\documentclass[letterpaper, 10pt, conference]{IEEEconf}
\IEEEoverridecommandlockouts 

\usepackage{amsmath,amssymb}
\usepackage{epsfig,subfigure,url}

\newcommand{\norm}[1]{\ensuremath{\left\| #1 \right\|}}
\newcommand{\bracket}[1]{\ensuremath{\left[ #1 \right]}}
\newcommand{\braces}[1]{\ensuremath{\left\{ #1 \right\}}}
\newcommand{\parenth}[1]{\ensuremath{\left( #1 \right)}}
\newcommand{\refeqn}[1]{(\ref{eqn:#1})}
\newcommand{\reffig}[1]{Fig. \ref{fig:#1}}
\newcommand{\tr}[1]{\mbox{tr}\ensuremath{\bracket{#1}}}

\newcommand{\SO}{\ensuremath{{SO(3)}}}

\renewcommand{\Re}{\ensuremath{\mathbb{R}}}
\renewcommand{\S}{\ensuremath{\mathbb{S}}}

\graphicspath{{./figs/}}

\title{\LARGE \bf
Propagation of Uncertainty in Rigid Body Attitude Flows}

\author{Taeyoung Lee\authorrefmark{1}\authorrefmark{2}, Nalin A. Chaturvedi\authorrefmark{2}, Amit Sanyal, Melvin Leok\authorrefmark{1}, and N. Harris McClamroch\authorrefmark{2}%
\thanks{Taeyoung Lee, Nalin A. Chaturvedi, and N. Harris McClamroch, Aerospace Engineering, University of Michigan, Ann Arbor, MI 48109 {\tt \{tylee,nalin,nhm\}@umich.edu}}%
\thanks{Amit Sanyal, Mechanical Engineering, University of Hawaii at Manoa, Honolulu, HI  96822  {\tt aksanyal@hawaii.edu}}%
\thanks{Melvin Leok, Mathematics, Purdue University, West Lafayette, IN 47907 {\tt mleok@math.purdue.edu}}%
\thanks{\textsuperscript{\footnotesize\ensuremath{*}}This research has been supported in part by NSF under grant DMS-0504747, and by a grant from the Rackham Graduate School, University of Michigan.}
\thanks{\textsuperscript{\footnotesize\ensuremath{\dagger}}This research has been supported in part by NSF under grant CMS-0555797.}
}

\begin{document}
\maketitle \thispagestyle{empty} \pagestyle{empty}

\begin{abstract}
Motivated by attitude control and attitude estimation problems for a rigid body, computational methods are proposed to propagate uncertainties in the angular velocity and the attitude.   The nonlinear attitude flow is determined by Euler-Poincar\'e equations that describe the rotational dynamics of the rigid body acting under the influence of an attitude dependent potential and by
a reconstruction equation that describes the kinematics expressed in terms of an orthogonal matrix representing the rigid body attitude.  Uncertainties in the angular velocity and attitude are described in terms of ellipsoidal sets that are propagated through this highly nonlinear attitude flow.  Computational methods are proposed, one method based on a local linearization of the attitude flow and two methods based on propagation of a small (unscented) sample selected from the initial uncertainty ellipsoid.   Each of these computational methods is constructed using the Lie group variational integrator algorithm, viewed as a discretization of the attitude flow dynamics.   Computational results are obtained that indicate (1) the strongly nonlinear attitude flow characteristics and (2) the limitations of each of these methods, and indeed any method, in providing effective global bounds on the nonlinear attitude flow.
\end{abstract}

\section{Introduction}

As an integrable system, the attitude dynamics of the free rigid body are reasonably well understood.    However, if there is an attitude dependent potential that influences the rigid body, then the dynamics can be very complex.    In this paper, such attitude dynamics are studied for a rigid body with an inertially fixed pivot acting under the influence of uniform and constant gravity; this model is subsequently referred to as the 3D pendulum~\cite{SheSanCha.CDC04}.

A rigid 3D pendulum is a rigid body supported by a fixed, frictionless pivot, acted on by gravitational forces. The supporting pivot allows three degrees of rotational freedom of the pendulum. Uniform, constant gravity is assumed. The terminology 3D pendulum refers to the fact that the pendulum is a rigid body with three spatial dimensions and the pendulum has three rotational degrees of freedom.

Two reference frames are introduced. An inertial reference frame has its origin at the pivot; the first two axes lie in the horizontal plane and the third axis is vertical in the direction of gravity. A reference frame fixed to the pendulum body is also introduced.  The origin of this body-fixed frame is also located at the pivot.  In the body fixed frame, the moment of inertia of the pendulum, $J\in\Re^{3\times 3}$, is constant. This moment of inertia can be computed from the traditional moment of inertia of a translated frame whose origin is located at the center of mass of the pendulum using the parallel axis theorem.  Since the origin of the body fixed frame is located at the pivot, principal axes with respect to this frame can be defined for which the moment of inertia is a diagonal matrix.  Note that the center of mass of the 3D pendulum may or may not lie on one of the principal axes defined in this way.

The dynamics of the 3D pendulum are given by the Euler equation that includes the moment due to gravity:
\begin{equation}
J \dot\Omega =J\Omega \times \Omega + mg \rho \times R^T e_3. \label{eq:Jw_dot} \end{equation}
The rotational kinematics equations are
\begin{equation} \dot{R} = R S(\Omega).\label{eq:R_dot} \end{equation}
Equations (\ref{eq:Jw_dot}) and (\ref{eq:R_dot}) define the full dynamics of a
rigid pendulum on the tangent bundle $TSO(3)$. The attitude of the pendulum is represented by a rotation matrix $R\in SO(3)$, and the angular velocity in the body fixed frame is denoted by $\Omega\in\Re^3$. 

In the above equations, the mass of the pendulum and the gravitational acceleration are denoted by $m$ and $g$, respectively, and $\rho\in\Re^3$ is the vector from the pivot to the mass center of the pendulum expressed in the body fixed frame. The vector $e_3 = [0,0,1]^T$ is the direction of gravity in the inertial frame, so that
$R^T e_3$ is the direction of gravity in the pendulum-fixed frame.
For a given vector $a\in\Re^3$, the skew-symmetric matrix $S(a)$ is defined as
\begin{equation}
S(a) = \begin{bmatrix} 0 & -a_3 & a_2 \\ a_3 &0& -a_1 \\ -a_2
&a_1& 0
\end{bmatrix}
\end{equation}
so that $S(a)b=a\times b$ for a $b\in\Re^3$.

There are two disjoint equilibria when the direction of gravity in the body fixed frame is collinear with the vector $\rho$. We define
\begin{align*}
    H=\braces{(R,\Omega)\in TSO(3)\,\big|\, R^Te_3=\frac{\rho}{\norm{\rho}},\Omega=0},\\
    I=\braces{(R,\Omega)\in TSO(3)\,\big|\, R^Te_3=-\frac{\rho}{\norm{\rho}},\Omega=0}
\end{align*}
as hanging equilibria and inverted equilibria, respectively.

The objective of this paper is to study the nonlinear attitude flow of the 3D pendulum dynamics by characterizing the propagation of uncertainty in the attitude and angular velocity.   Simple analytical bounds can be obtained for small, localized uncertainty over short time periods, but such bounds provide no insight into the global attitude dynamics.   Uncertainty propagation can also be studied, within a probabilistic framework, using the Liouville partial differential equation, a special case of the Fokker-Planck equation, but this leads to significant solution difficulties.  Our approach in this paper is to make use of computational tools, based on the recently introduced Lie group variational integrator for the 3D pendulum and on a framework for propagating suitably defined uncertainty ellipsoids~\cite{CCA05}.  We show that these computational tools are useful in some cases, but they have important limitations that arise from the complex attitude dynamics of the 3D pendulum.

This line of research was motivated by our prior work on attitude control and estimation~\cite{CDC06.est,SCL06}.   The results in this paper demonstrate the complex global dynamics that can occur for an uncontrolled 3D pendulum.    In this way, the challenges of attitude control and estimation are made clear, at least in the case that global results are desired.

\section{Problem formulation}
An uncertainty ellipsoid on $T\SO$ is defined by
\begin{align*}
    \mathcal{E}(\hat R,\hat\Omega,P)=\braces{(R,\Omega)\in T\SO\,\big|\,
    x^TP^{-1}x\leq 1},
\end{align*}
where $x=[\zeta;\delta\Omega]\in\Re^6$ and $\zeta=\log (\hat R^T
R),\delta\Omega=\Omega-\hat\Omega\in\Re^3$. The center of the ellipsoid is given by the rotation matrix $\hat R$ and the angular velocity $\hat\Omega$; the matrix $P\in\Re^{6\times 6}$ is the
uncertainty matrix that characterizes the size of the ellipsoid.   
In particular, if the initial attitude and angular velocity are known to lie in an initial uncertainty ellipsoid, we seek computational methods to propagate the uncertainty so 
that we obtain another uncertainty ellipsoid at a later time in which the resulting attitude and angular velocity of the 3D pendulum are expected to lie.

\section{Uncertainty propagation}
We study three methods to propagate the uncertainty set over the time interval $[0, T]$.  All three methods make use of the Lie group variational
integrator to propagate the angular velocity and the attitude~\cite{CCA05}. The Lie group variational integrator is described by the following discrete update equations,
\begin{gather}
h S(J\Omega_k+\frac{h}{2} M_k) = F_k J_d - J_dF_k^T,\label{eqn:findf0}\\
R_{k+1} = R_k F_k,\label{eqn:updateR0}\\
J\Omega_{k+1} = F_k^T J\Omega_k +\frac{h}{2} F_k^T M_k
+\frac{h}{2}M_{k+1},\label{eqn:updatew0}
\end{gather}
where the subscript $k$ denotes the $k$th discrete variable for a
fixed integration step size $h\in\Re$, and $F_k\in\SO$ is the
relative attitude between two adjacent integration steps. The nonstandard moment inertia matrix is given by $J_d=\frac{1}{2}\tr{J}I_{3\times 3}-J\in\Re^{3\times 3}$, and the moment due to the gravity is denoted by $M_k=mg\rho\times R_k^T e_3\in\Re^3$. For a
given $(R_k,\Omega_k)$, \refeqn{findf0} is solved to
find $F_k\in\SO$. Then $(R_{k+1},\Omega_{k+1})$ is obtained by
\refeqn{updateR0} and \refeqn{updatew0}. This yields a map
$(R_k,\Omega_k)\mapsto(R_{k+1},\Omega_{k+1})$ and this process is
repeated.

\subsection{Uncertainty propagation using linearization on $[0, T]$}  In this method, the uncertainty ellipsoid is propagated by updating both the center of the ellipsoid and the the uncertainty matrix $P$.
The center of the uncertainty ellipsoid, denoted by $\hat R_k,\hat\Omega_k$, is propagated according to the Lie group variational integrator, initialized by the center of the initial uncertainty ellipsoid.
We then assume that the uncertainty ellipsoid is sufficiently small that the flow of points within the uncertainty ellipsoid is well-approximated by the linearized flow of the Lie group variational integrator about this center solution, which we denote by
\begin{align*}
x_{k+1} & = A_k x_k,
\end{align*}
where $x_k=[\zeta_k;\delta\Omega_k]\in\Re^6$ and the matrix $A_k\in\Re^{6\times 6}$ depends on $\hat R_k,\hat \Omega_k$. The expression for $A_k$ can be found in~\cite{ACC06}. Using this linearization, the uncertainty matrix is propagated according to
\begin{align}
P_{k+1}& = A_k P_k A_k^T\label{eqn:Pkpf}
\end{align}
with initial condition given by the initial uncertainty matrix.  In this way, the uncertainty ellipsoid is propagated and determined at time $T$.

\subsection{Uncertainty propagation using unscented method on $[0,T]$} In this method, we compute the propagated uncertainty ellipsoid by using the Lie group variational integrator to propagate a small sample of points selected from the initial uncertainty ellipsoid.   Following the conventional wisdom, we choose the 12 points corresponding to the intersection of the boundary of the initial uncertainty ellipsoid
and its principal axes. This choice is informed by the fact that if the initial uncertainty ellipsoid is propagated by a linear flow, it will remain an ellipsoid, and it will coincide with the minimal volume ellipsoid containing the the 12 intersection points propagated by the same linear flow.

More explicitly, suppose that the initial uncertainty
ellipsoid is given by
\begin{align*}
    \mathcal{E}(R_0,\Omega_0,P_0).
\end{align*}
Let $\lambda^i\in\Re$ and $\phi^i=[\phi_R^i;\phi_\Omega^i]\in\Re^6$ be the $i$-th
eigenvalue and eigenvector, respectively, of the uncertainty matrix $P_0\in\Re^{6\times 6}$ for $i\in\braces{1,2\ldots,6}$. The intersection of
the corresponding principal axis and the ellipsoid boundary is $\sqrt{\lambda_i}\phi_i$.
Then, the 12 intersection points are given by,
\begin{align*}
    \braces{\parenth{R_0\exp
    (\pm\sqrt{\lambda^i}\hat\phi^i_R),\Omega_0\pm\sqrt{\lambda^i}\phi^i_\Omega}
    }
    \quad i\in\braces{1,2,\ldots,6}.
\end{align*}
We propagate these 12 initial conditions using the Lie group variational integrator to determine the 12 values of the attitude and angular velocity at time $T$.    We then determine a minimal volume ellipsoid that contains all 12 values of the attitude and angular velocity at time $T$~\cite{BoyVan.BK04}.  From this computed ellipsoid, the center attitude, the center angular velocity, and the uncertainty matrix at time $T$ can be computed.

\subsection{Uncertainty propagation using unscented method with re-sampling on $[0,T]$}
The unscented method propagates the same sampled points obtained from the \textit{initial}
uncertainty ellipsoid throughout  the entire time period $[0,T]$.   In this modification, we partition the time period
$[0,T]$ into subintervals and we use the Lie group variational integrator to propagate the 12 sample points throughout each subinterval.  At the end of each
subinterval, we determine a minimal volume ellipsoid that contains all 12 values of the attitude and angular velocity at the end of that subinterval.    We then select a new set of sample
points located on the principal axes of this new uncertainty ellipsoid. On the first subinterval, the 12 initial values of the attitude and angular velocity are obtained from the initial uncertainty ellipsoid.
This method differs from the
previous method in that resampling of the propagated 12 solutions is carried out at the end of each subinterval.

It is worth remarking that if the flow is linear, the re-sampling technique will yield the same estimates as the usual unscented method. However, if the flow is nonlinear, and the re-sampling is sufficiently frequent, then the re-sampling technique will tend to yield a more conservative estimate. Indeed, the difference between the two estimates gives an indication of how nonlinear the flow is.

\section{Numerical examples}
We apply these methods to the attitude dynamics of the 3D pendulum.
The pendulum body is chosen as an elliptic cylinder, and its properties are given by
\begin{align*}
    J=\mathrm{diag}[0.13,0.28,0.17]\,\mathrm{kg\cdot m^2},\; m=1\,\mathrm{kg},\;
    \rho=0.3e_3\,\mathrm{m}.
\end{align*}
We consider two initial conditions; the first initial condition results in a near-oscillatory attitude flow while the second initial condition results in a highly irregular attitude flow.

\begin{figure}
\centerline{
    \subfigure[Reduced attitude $\Gamma$]{
    \includegraphics[width=0.42\columnwidth]{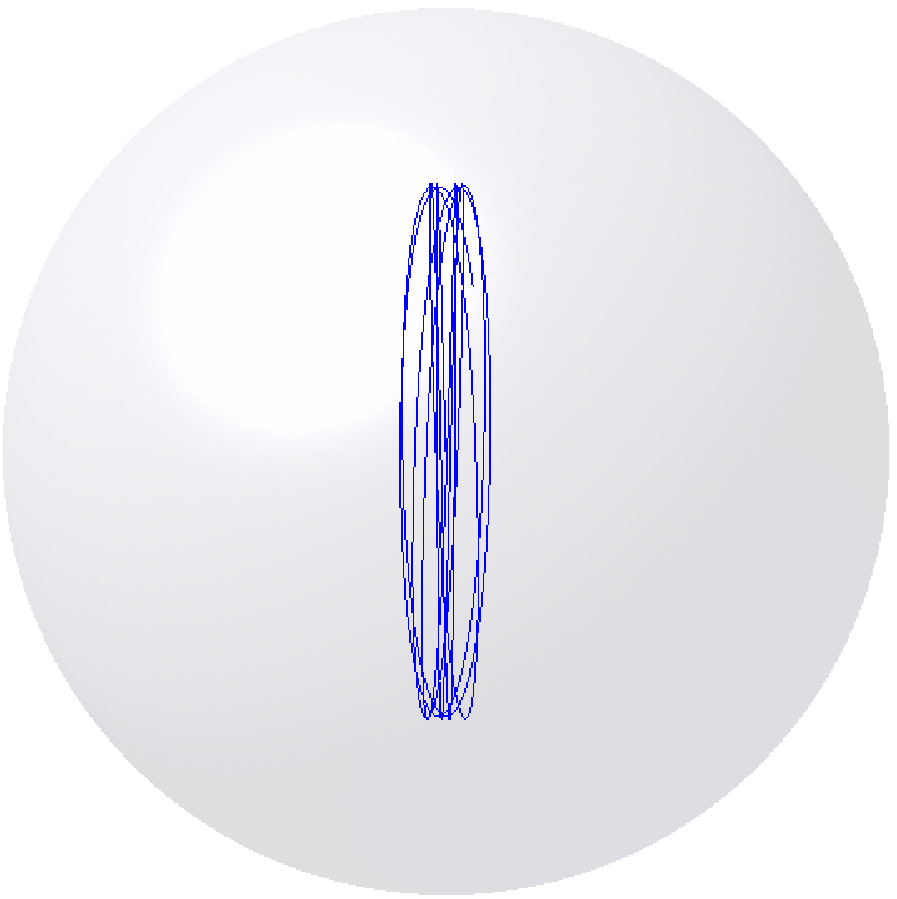}}
    \subfigure[Angular velocity $\Omega$]{
    \includegraphics[width=0.55\columnwidth]{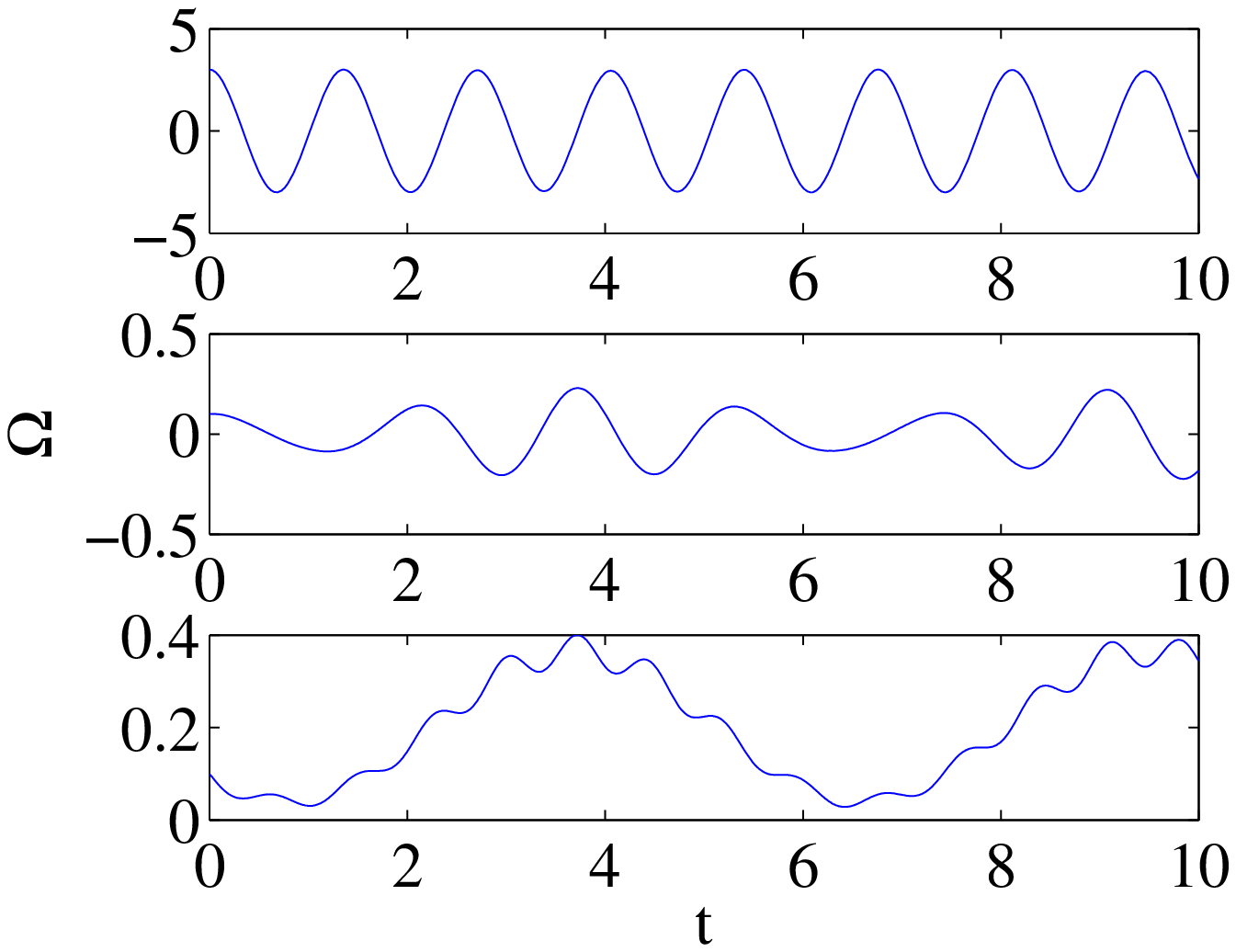}}
    }
\caption{Reduced attitude and angular velocity for oscillatory attitude
flow (The center of the sphere is the hanging equilibrium when $R^Te_3=e_3$.)}\label{fig:oscGamw}
\end{figure}

\paragraph*{Oscillatory attitude flow}
The initial uncertainty ellipsoid is defined by its initial center attitude and angular velocity and the initial uncertainty matrix, which are given by
\begin{gather*}
    R_0=I_{3\times 3},\quad \Omega_0=[3.0,0.1,0.1]\,\mathrm{rad/s},\\
    P_0=\mathrm{diag}\bracket{(5\frac{\pi}{180})^2[1,1,1],\,0.01^2[1,1,1]}.
\end{gather*}
It is convenient for graphical purposes to consider the reduced attitude $R^T e_3$ on the two-sphere $\S^2$ and to plot the induced attitude flow on $\S^2$. This is due to the fact that the 3D pendulum has a symmetry of $\S^1$ given by a rotation about the vertical axis, and the configuration space can be reduced to a quotient space $\SO/\S^1\simeq\S^2$ using this symmetry~\cite{SIADS07}. The reduced attitude denotes the direction of gravity in the body fixed frame. Plots of the reduced attitude and the angular velocity responses, corresponding to the initial conditions, are
shown in \reffig{oscGamw} for a time period of $10$ seconds.

We apply the three methods to propagate the initial uncertainty through this oscillatory
attitude flow for $10$ seconds. The integration step size is $h=0.005$, corresponding to 2000 time steps over the length of the simulation. For the unscented
method, we find the minimal volume covering ellipsoid every $0.1$ seconds. For the
unscented method with re-sampling, we find the minimal volume covering ellipsoid every
$0.1$ seconds, and we choose new set of sample points every $2.0$ seconds. To provide a baseline for comparing the performance of each method, we choose sample points in the interior of the initial uncertainty ellipsoid. In particular, we choose 144 points on the level
set of $x^TP^{-1}_0 x=0.8$, and we numerically integrate each of them.

\begin{figure}
\centerline{
    \subfigure[Magnitude of uncertainty matrix]{
    \includegraphics[width=0.50\columnwidth]{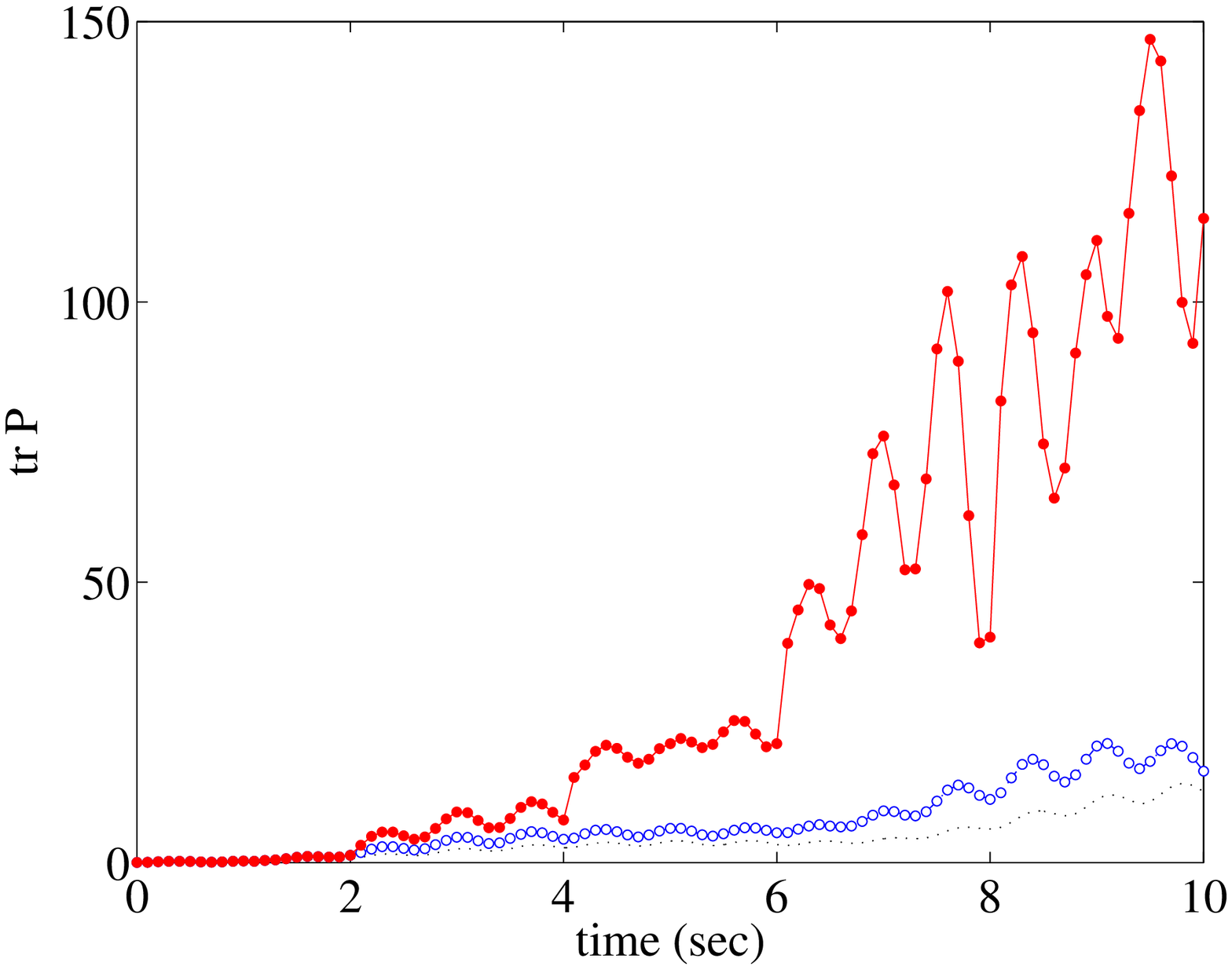}}
    \subfigure[Percentage of sample points contained in the computed uncertainty (\%)]{
    \includegraphics[width=0.50\columnwidth]{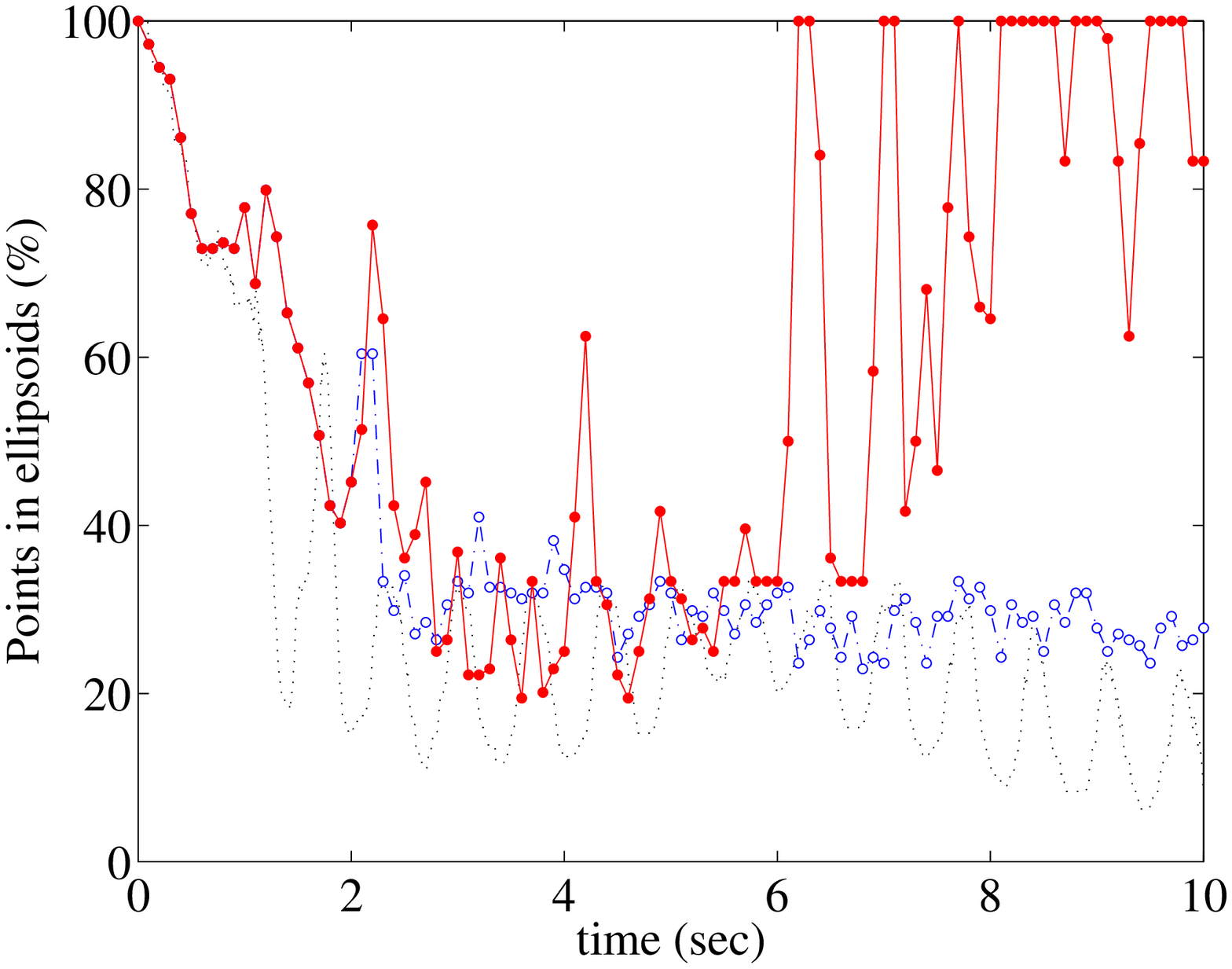}}
    }
\caption{Uncertainty propagation through oscillatory attitude flow (Linearization: dotted, Unscented: blue, Unscented with re-sampling: red)}\label{fig:oscun}
\end{figure}

\begin{figure*}
\footnotesize\selectfont
\centerline{
    \subfigure{\includegraphics[width=0.16\textwidth]{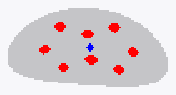}}
    \hspace*{0.02\textwidth}
    \subfigure{\includegraphics[width=0.16\textwidth]{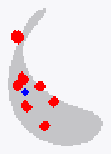}}
    \hspace*{0.02\textwidth}
    \subfigure{\includegraphics[width=0.16\textwidth]{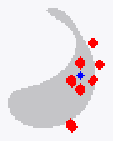}}
    \hspace*{0.02\textwidth}
    \subfigure{\includegraphics[width=0.16\textwidth]{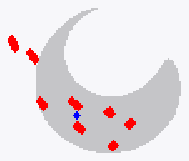}}
    \hspace*{0.02\textwidth}
    \subfigure{\includegraphics[width=0.16\textwidth]{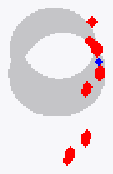}}
}
\centerline{Linearization}\vspace*{0.2cm}
\centerline{
    \subfigure{\includegraphics[width=0.16\textwidth]{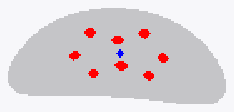}}
    \hspace*{0.02\textwidth}
    \subfigure{\includegraphics[width=0.16\textwidth]{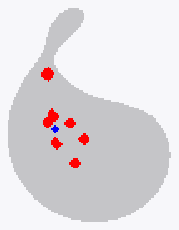}}
    \hspace*{0.02\textwidth}
    \subfigure{\includegraphics[width=0.16\textwidth]{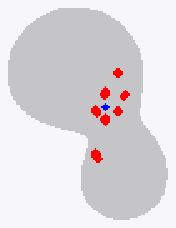}}
    \hspace*{0.02\textwidth}
    \subfigure{\includegraphics[width=0.16\textwidth]{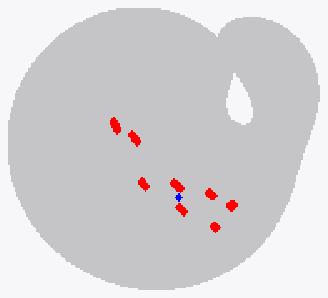}}
    \hspace*{0.02\textwidth}
    \subfigure{\includegraphics[width=0.16\textwidth]{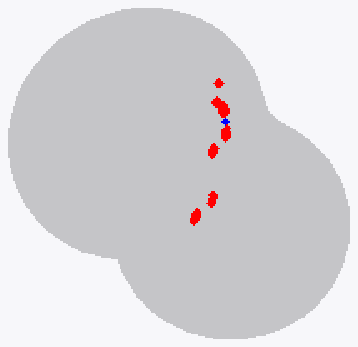}}
}
\centerline{Unscented method}\vspace*{0.2cm}
\centerline{
    \subfigure[$t=2.0$]{\includegraphics[width=0.16\textwidth]{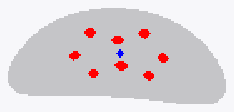}}
    \hspace*{0.02\textwidth}
    \subfigure[$t=4.0$]{\includegraphics[width=0.16\textwidth]{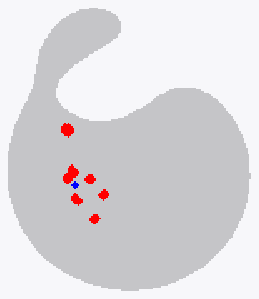}}
    \hspace*{0.02\textwidth}
    \subfigure[$t=6.0$]{\includegraphics[width=0.16\textwidth]{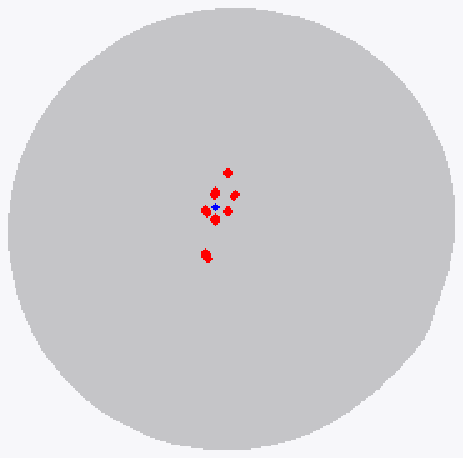}}
    \hspace*{0.02\textwidth}
    \subfigure[$t=8.0$]{\includegraphics[width=0.16\textwidth]{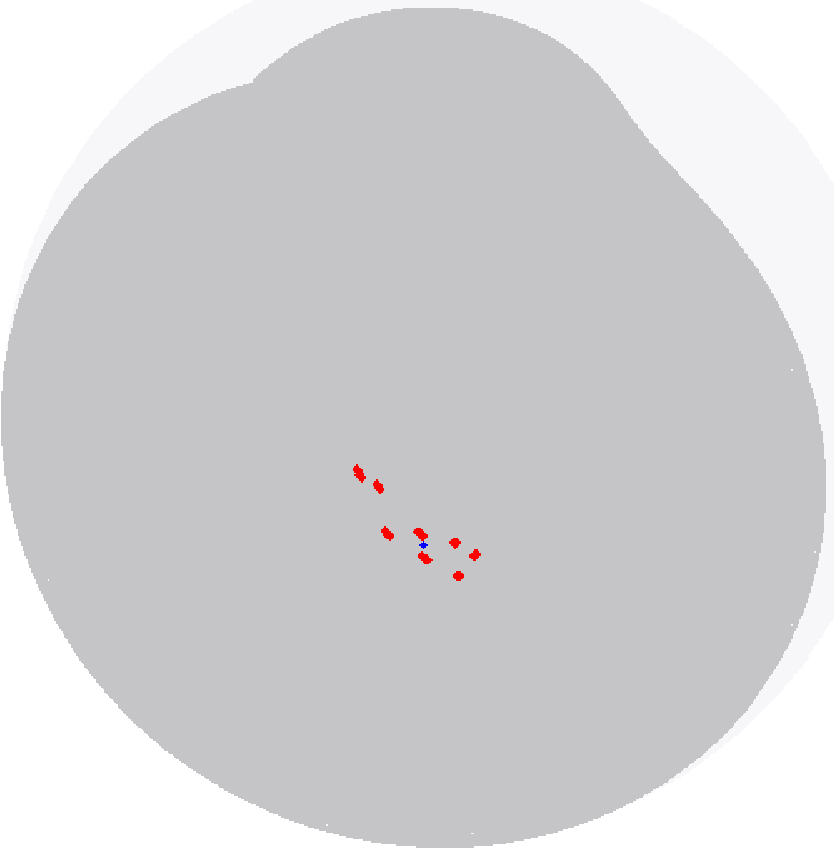}}
    \hspace*{0.02\textwidth}
    \subfigure[$t=10.0$]{\includegraphics[width=0.16\textwidth]{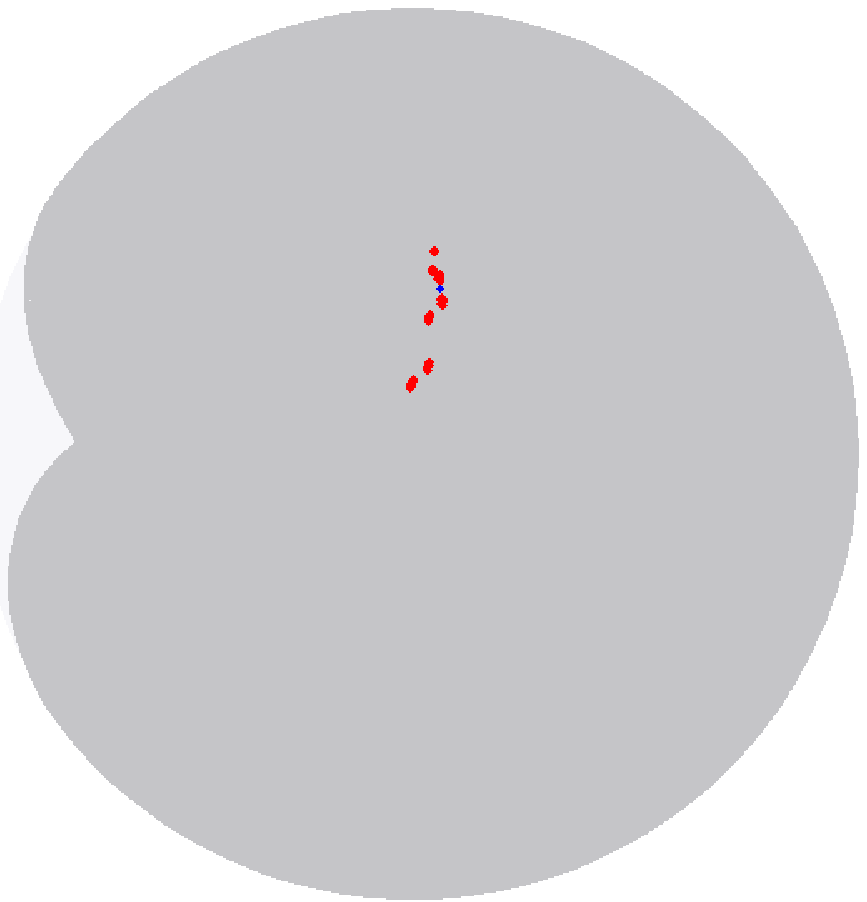}}
}
\centerline{Unscented with re-sampling}
\caption{Uncertainty projected onto the reduced attitude on $\S^2$ for oscillatory attitude flow (The center of the sphere is the hanging equilibrium when $R^Te_3=e_3$.)}\label{fig:oscsph}
\end{figure*}

The magnitudes of the uncertainty matrix and the percentage of sample points (out of 144) contained in the computed uncertainty ellipsoids are shown in \reffig{oscun}. The
propagation of uncertainties in the reduced attitude on $\S^2$ are shown in
\reffig{oscsph}.

For a short time period, the properties of all methods are similar. The uncertainty
ellipsoid computed by the unscented method is larger than for the linearization method, which is a reflection of the nonlinear nature of the dynamics. The re-sampling in the unscented method with re-sampling has the
effect of enlarging the uncertainty ellipsoids. Thus, the magnitude of the uncertainty
ellipsoid increases rapidly, and the uncertainty ellipsoid for the re-sampling method contains a greater proportion of the sample points after 6 seconds.

The mean precentages of sample points contained in the computed uncertainty ellipsoids are
28.9\%, 39.8\%, and 59.8\%, respectively. This suggests that even for this oscillatory
attitude flow, the nonlinear effects are so strong that it is difficult to accurately propagate the
uncertainty.

The computation times are 3.24, 45.18, and 45.84 seconds, respectively.


\begin{figure}
\centerline{
    \subfigure[Reduced attitude $\Gamma$]{
    \includegraphics[width=0.42\columnwidth]{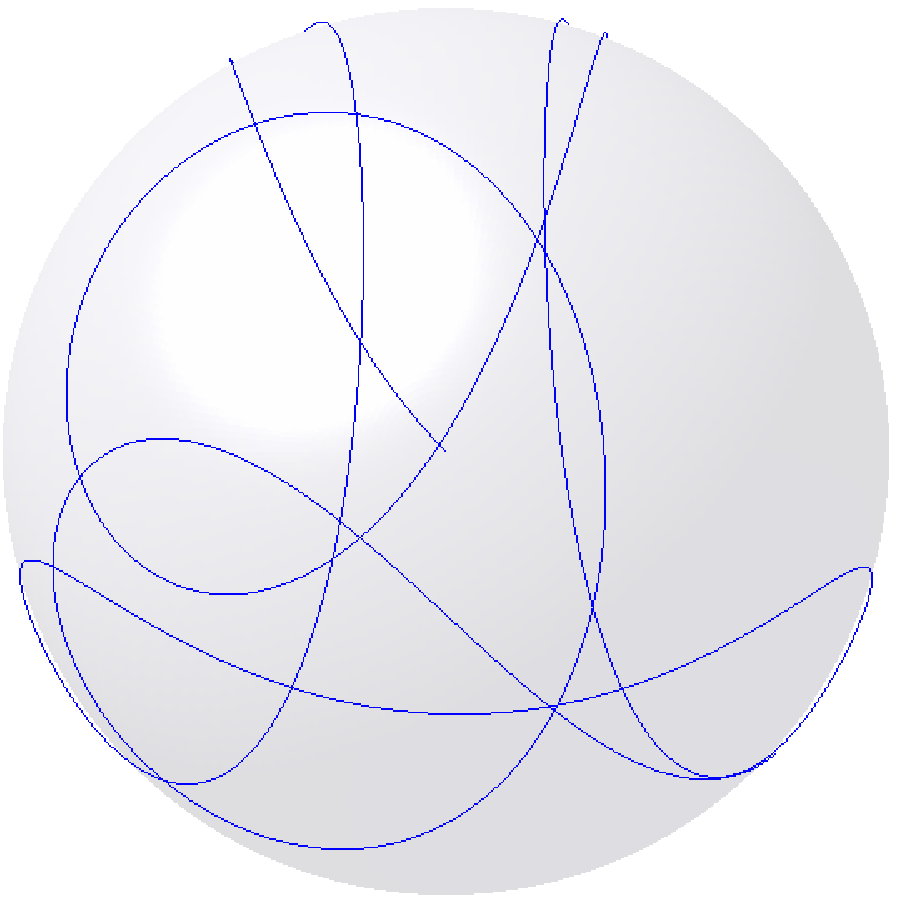}}
    \subfigure[Angular velocity $\Omega$]{
    \includegraphics[width=0.55\columnwidth]{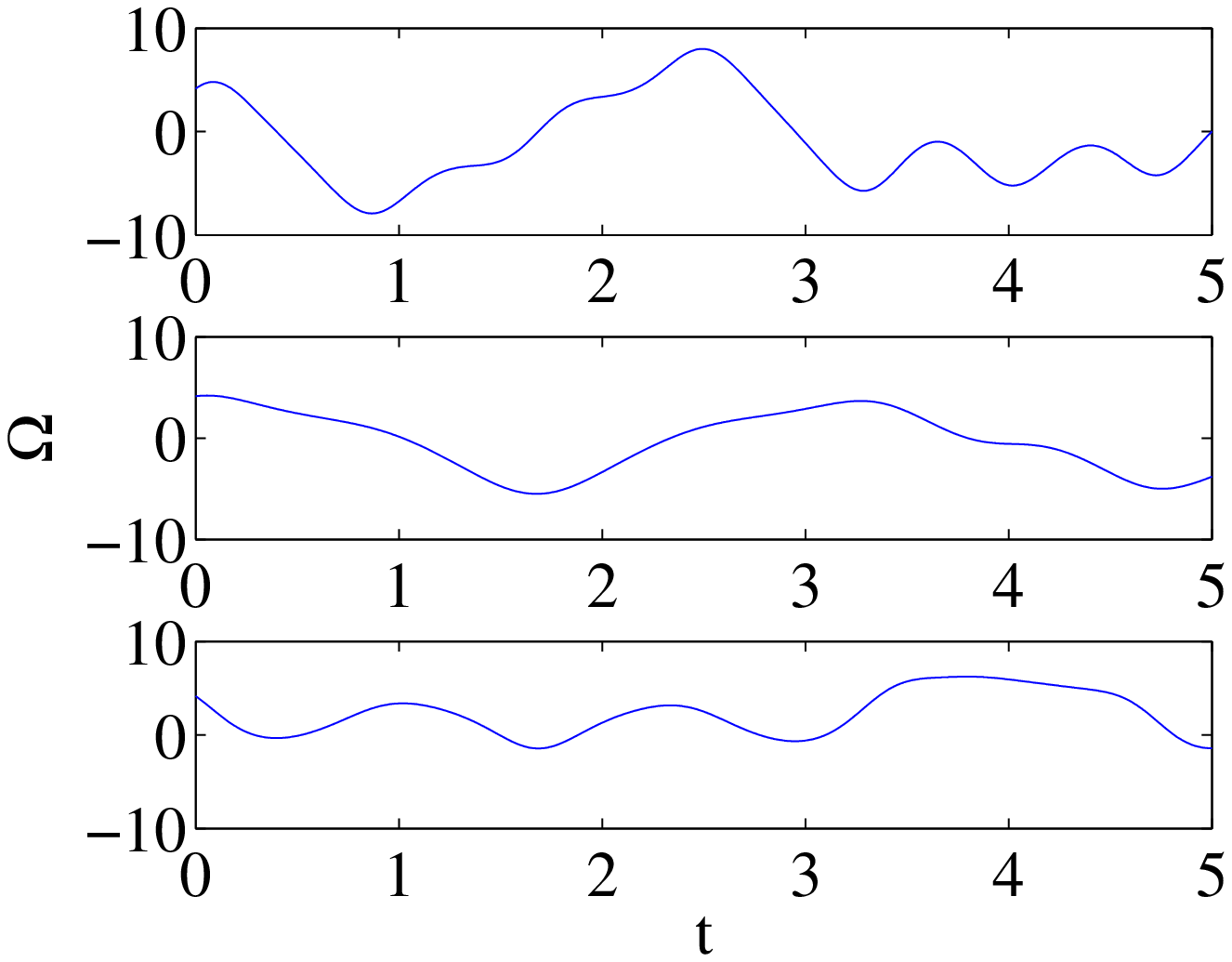}}
    }
\caption{Reduced attitude and angular velocity for irregular attitude
flow (The center of the sphere is the hanging equilibrium when $R^Te_3=e_3$.)}\label{fig:chaGamw}
\end{figure}

\paragraph*{Irregular attitude flow}
The initial uncertainty ellipsoid is defined by its initial attitude and angular velocity and the initial uncertainty matrix, which are given by
\begin{gather*}
    R_0=I_{3\times 3},\quad \Omega_0=[4.14,4.14,4.14]\,\mathrm{rad/s},\\
    P_0=\mathrm{diag}\bracket{(5\frac{\pi}{180})^2[1,1,1],\,0.01^2[1,1,1]}.
\end{gather*}
The reduced attitude on $\S^2$ and the angular velocity responses for $10$ seconds corresponding to the center initial conditions are
shown in \reffig{chaGamw}.

\begin{figure}
\centerline{
    \subfigure[Magnitude of uncertainty matrix]{
    \includegraphics[width=0.50\columnwidth]{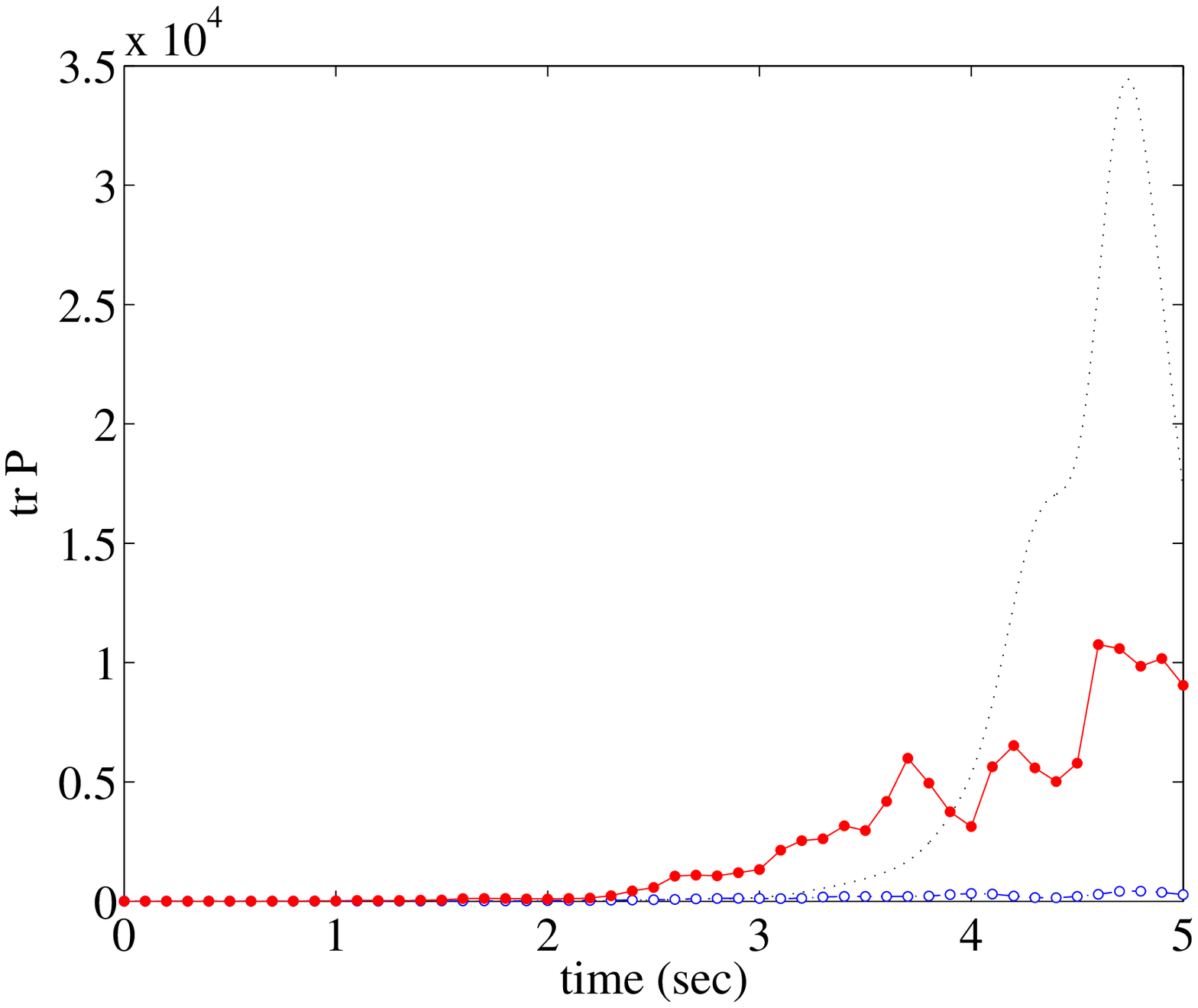}}
    \subfigure[Percentage of sample points contained in the computed uncertainty (\%)]{
    \includegraphics[width=0.50\columnwidth]{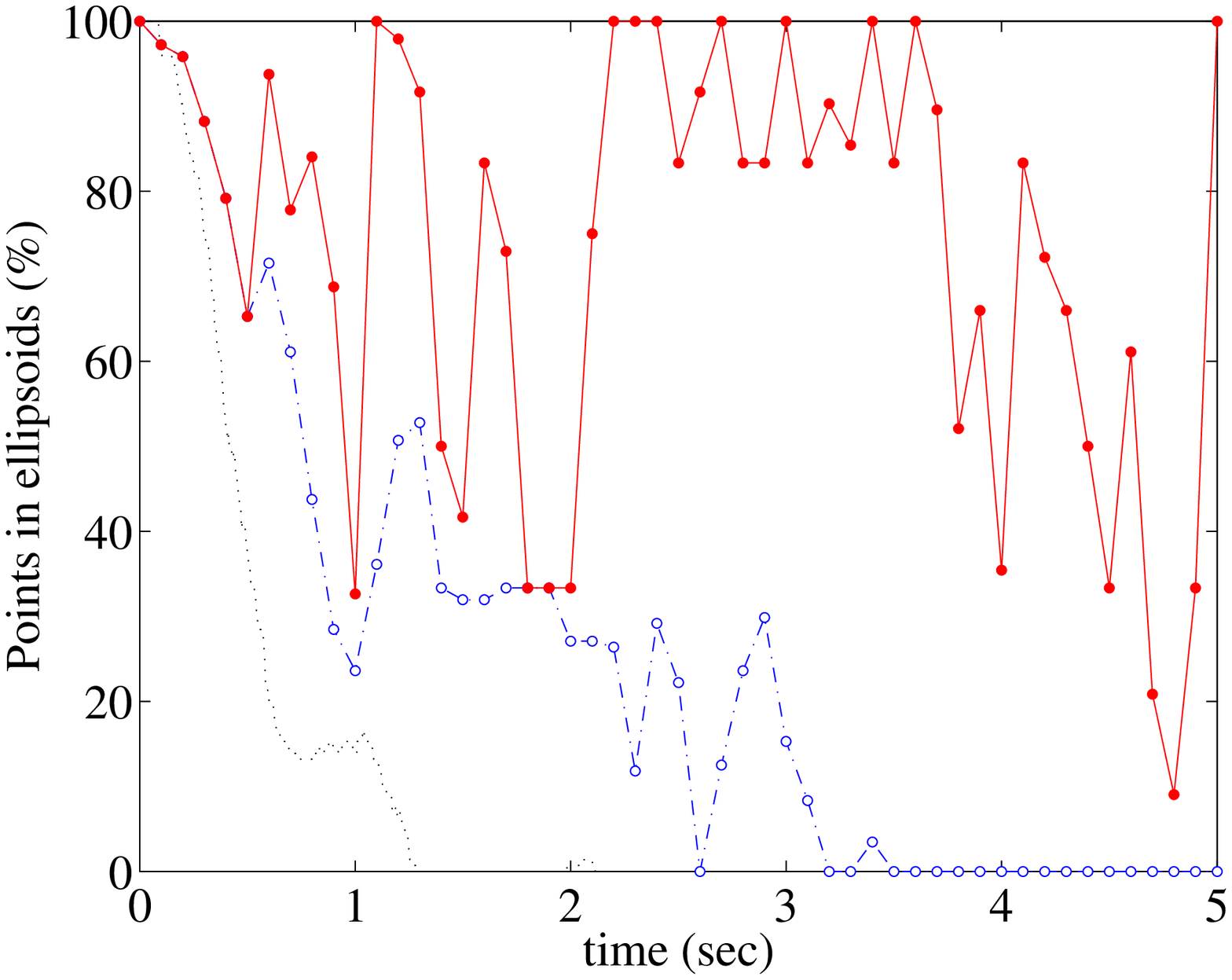}}
    }
\caption{Uncertainty propagation through irregular attitude flow (Linearization: dotted, Unscented: blue, Unscented with re-sampling: red)}\label{fig:chaun}
\end{figure}

We apply the three methods to propagate the initial uncertainty through this irregular
attitude flow for $5$ seconds. The integration step size is $h=0.002$, corresponding to 2500 time steps over the length of the simulation. For the unscented
method, we find the minimal volume covering ellipsoid every 0.1 seconds. For the
unscented method with re-sampling, we find the minimal volume covering ellipsoid every
$0.1$ seconds, and we choose new set of sample points every $0.5$ seconds. To compare the
properties of each method, we choose 144 initial sample points on the level set of
$x^TP^{-1}_0 x=0.8$, and we numerically integrate them.

\begin{figure*}
\footnotesize\selectfont
\centerline{
    \subfigure{\includegraphics[width=0.16\textwidth]{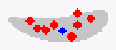}}
    \hspace*{0.02\textwidth}
    \subfigure{\includegraphics[width=0.16\textwidth]{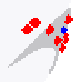}}
    \hspace*{0.02\textwidth}
    \subfigure{\includegraphics[width=0.16\textwidth]{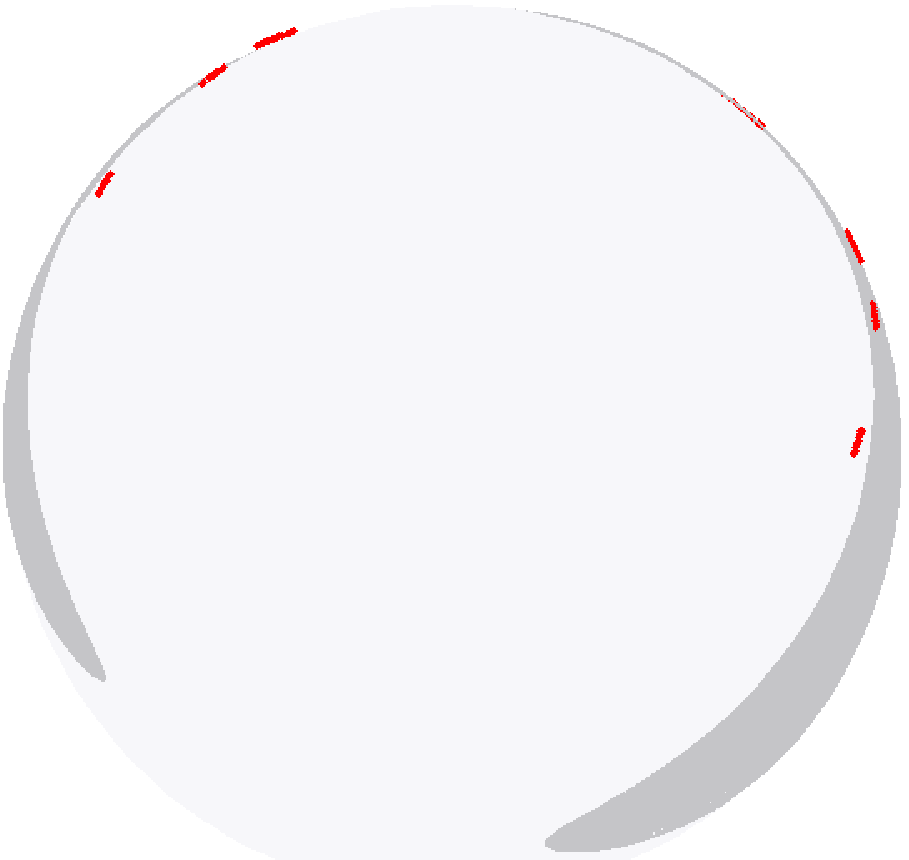}}
    \hspace*{0.02\textwidth}
    \subfigure{\includegraphics[width=0.16\textwidth]{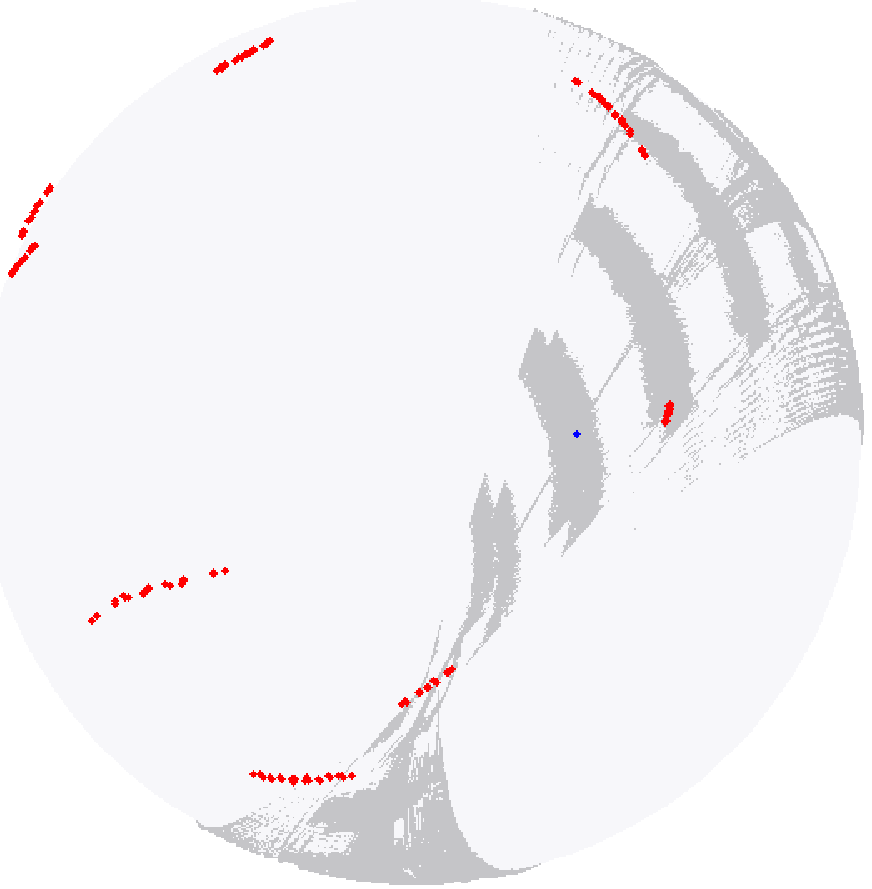}}
    \hspace*{0.02\textwidth}
    \subfigure{\includegraphics[width=0.16\textwidth]{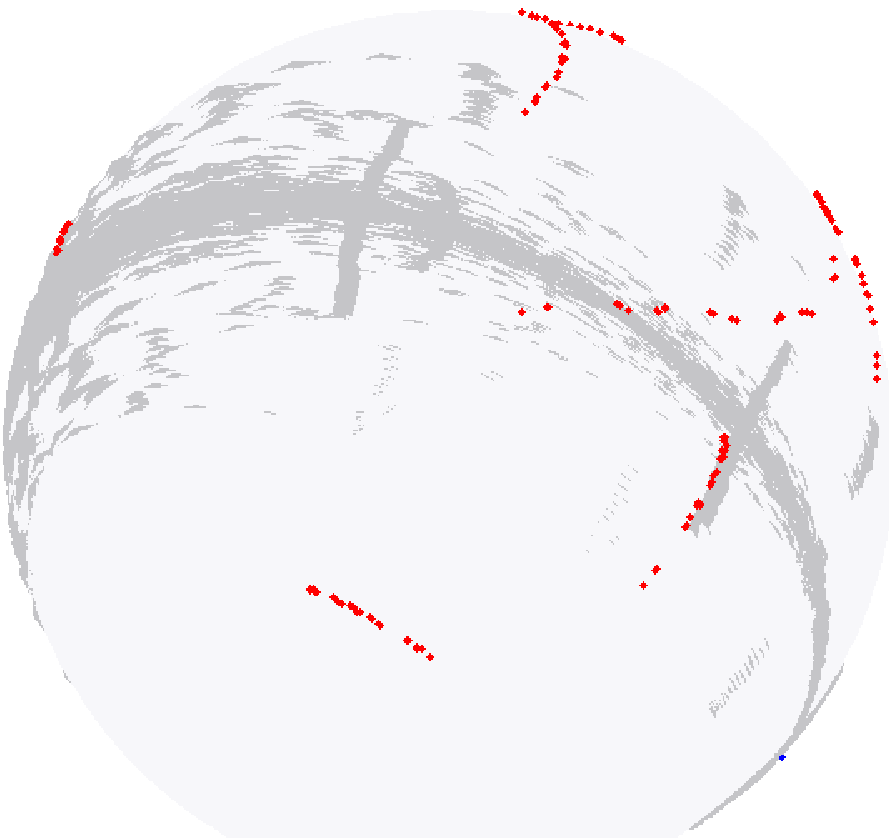}}
}
\centerline{Linearization}\vspace*{0.2cm}
\centerline{
    \subfigure{\includegraphics[width=0.16\textwidth]{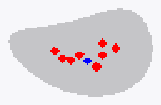}}
    \hspace*{0.02\textwidth}
    \subfigure{\includegraphics[width=0.16\textwidth]{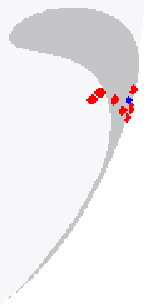}}
    \hspace*{0.02\textwidth}
    \subfigure{\includegraphics[width=0.16\textwidth]{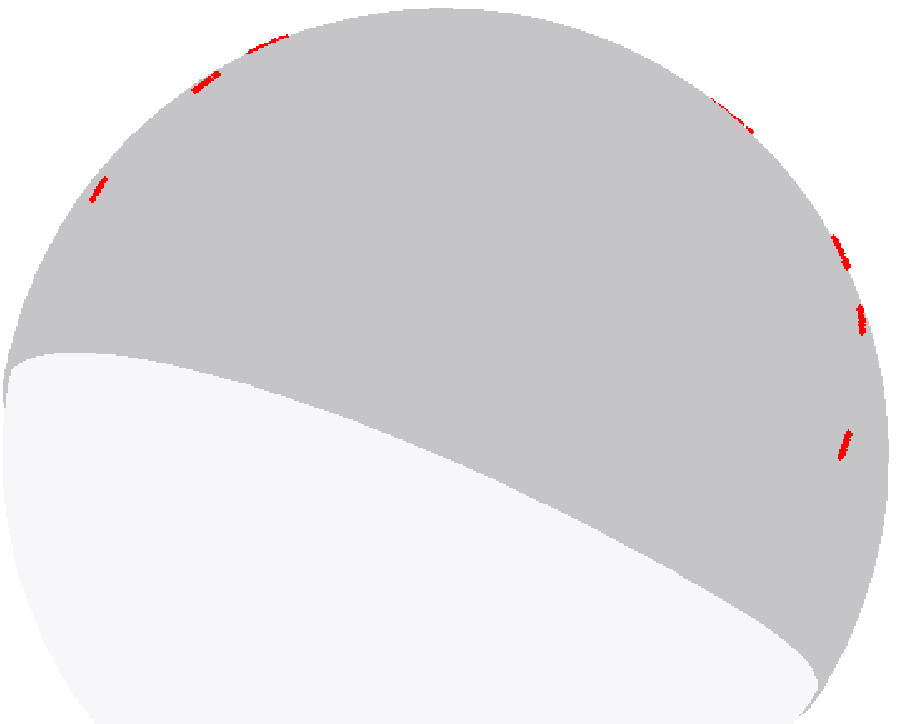}}
    \hspace*{0.02\textwidth}
    \subfigure{\includegraphics[width=0.16\textwidth]{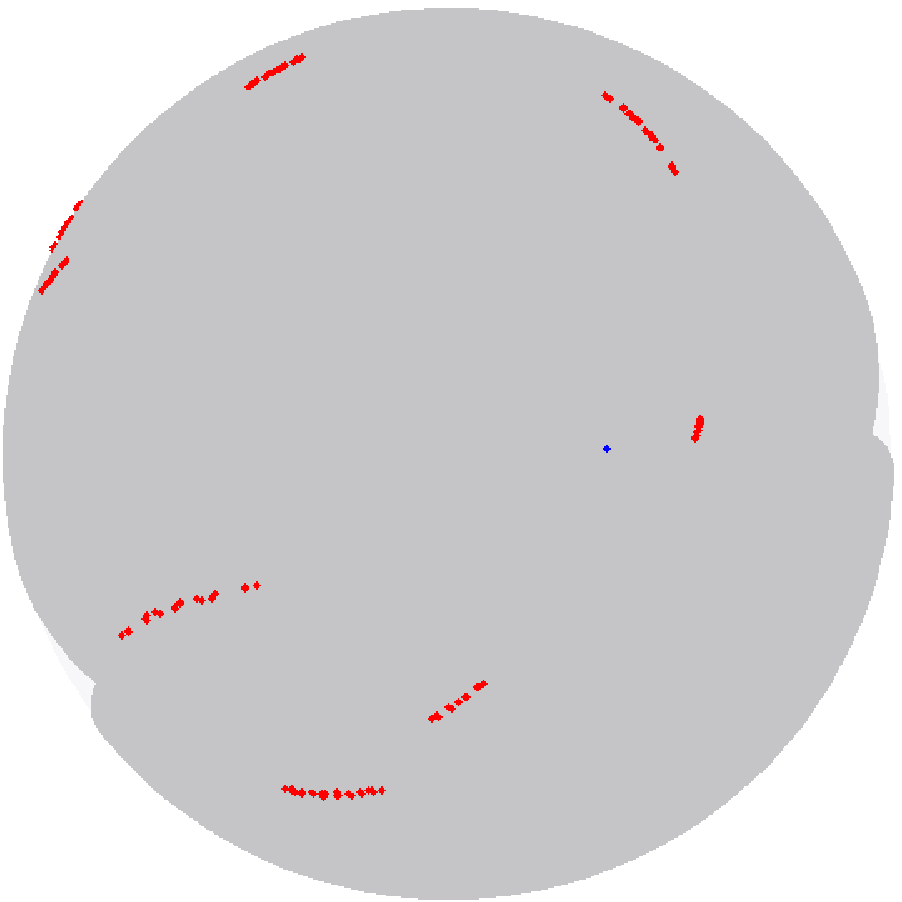}}
    \hspace*{0.02\textwidth}
    \subfigure{\includegraphics[width=0.16\textwidth]{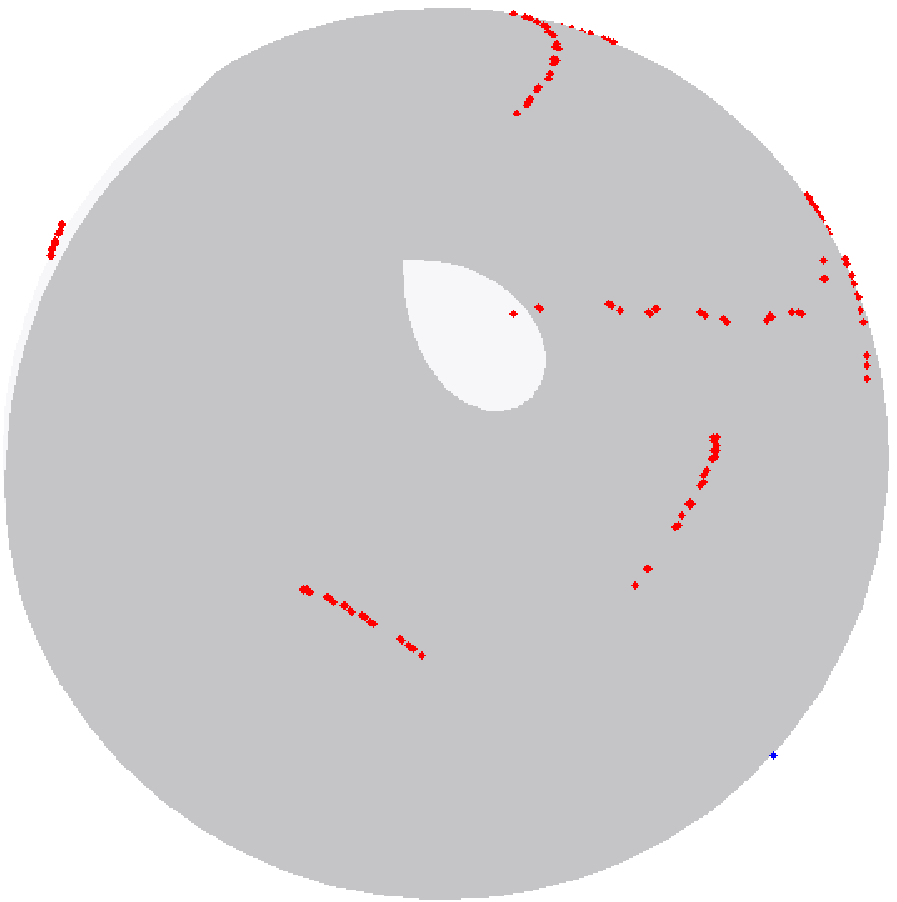}}
}
\centerline{Unscented method}\vspace*{0.2cm}
\centerline{
    \subfigure[$t=1.0$]{
    \includegraphics[width=0.16\textwidth]{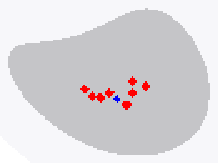}}
    \hspace*{0.02\textwidth}
    \subfigure[$t=2.0$]{
    \includegraphics[width=0.16\textwidth]{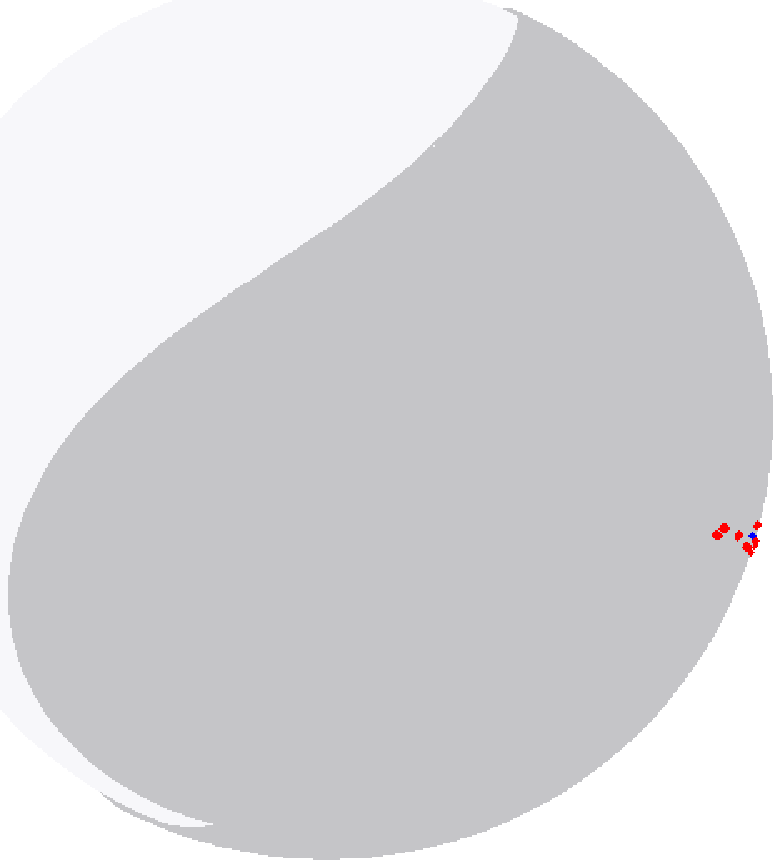}}
    \hspace*{0.02\textwidth}
    \subfigure[$t=3.0$]{
    \includegraphics[width=0.16\textwidth]{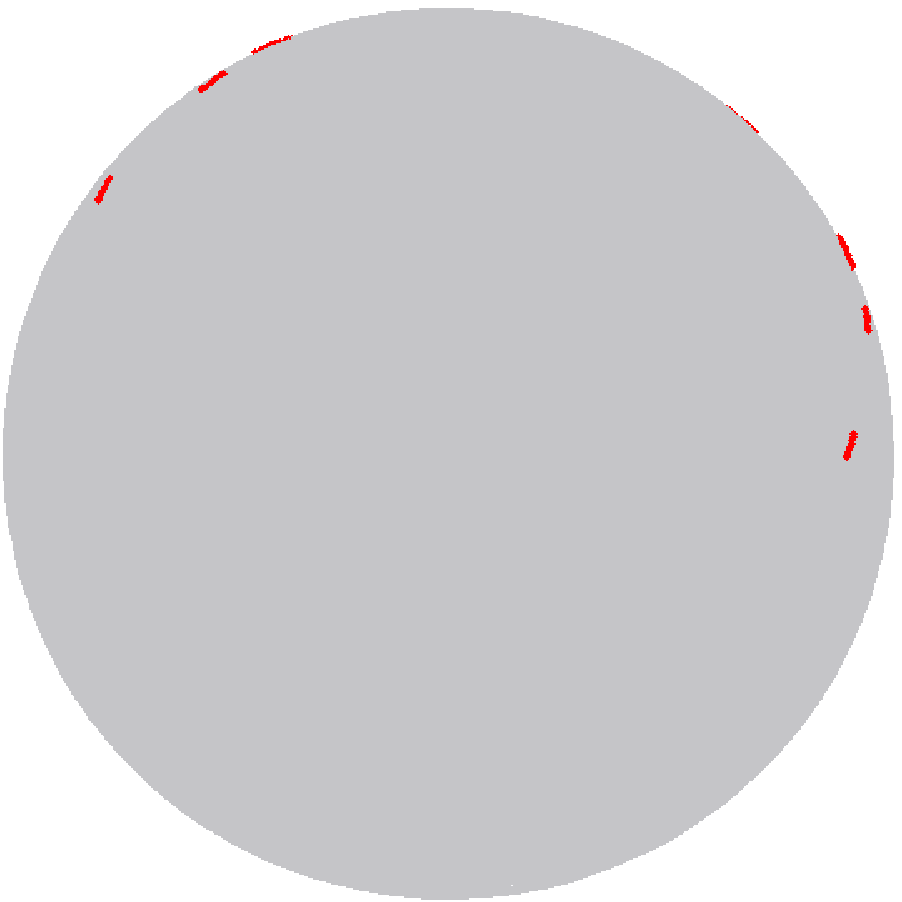}}
    \hspace*{0.02\textwidth}
    \subfigure[$t=4.0$]{
    \includegraphics[width=0.16\textwidth]{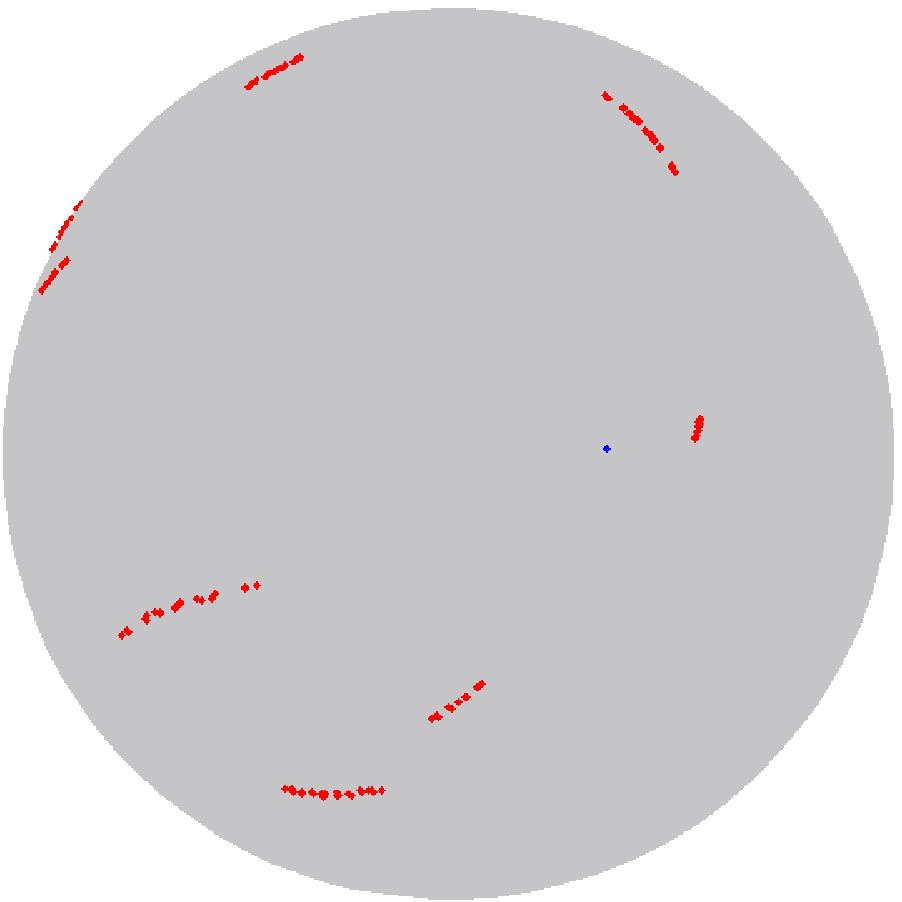}}
    \hspace*{0.02\textwidth}
    \subfigure[$t=5.0$]{
    \includegraphics[width=0.16\textwidth]{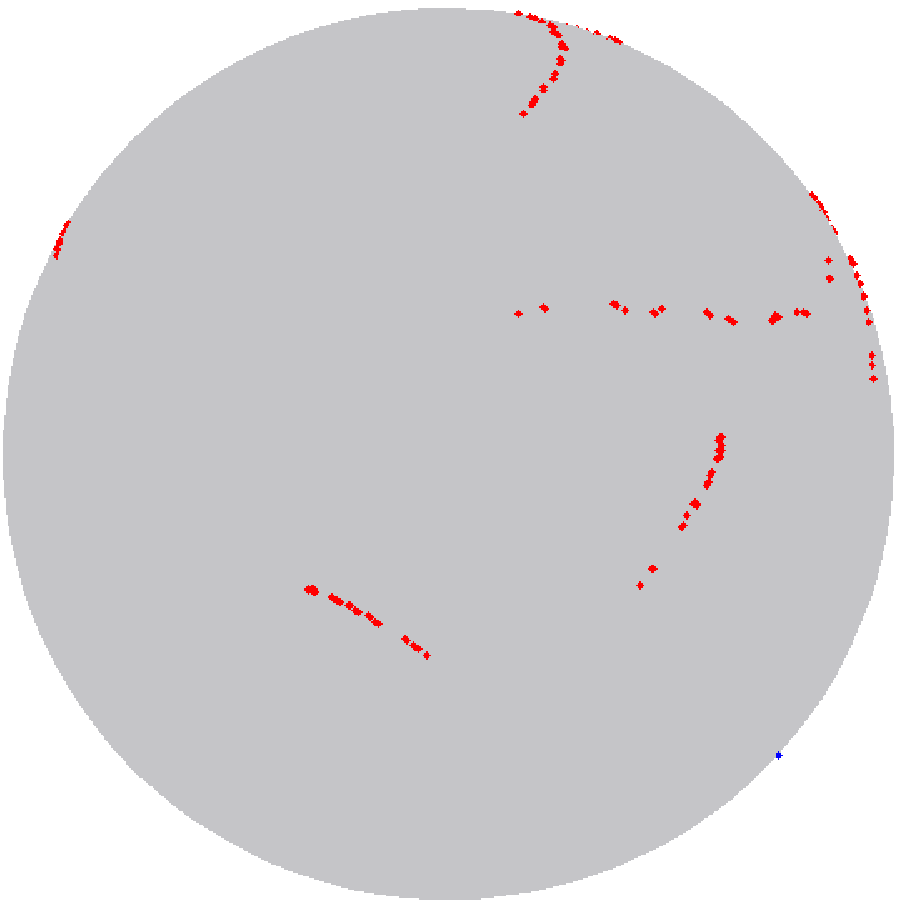}}
}
\centerline{Unscented with re-sampling}
\caption{Uncertainty projected onto the reduced attitude on $\S^2$ for irregular attitude flow (The center of the sphere is the hanging equilibrium when $R^Te_3=e_3$.)}\label{fig:chasph}
\end{figure*}

The magnitudes of the uncertainty matrix and the percentage of sample points (out of 144) contained in the computed uncertainty ellipsoid are shown in \reffig{chaun}. The
propagation of uncertainties in the reduced attitude on $\S^2$ are shown in
\reffig{chasph}.

For a short time period, the properties of all methods are similar. The uncertainty
ellipsoid computed by the linearization method grows rapidly, but it encloses few points after 1.5
second. The unscented method encloses more points with smaller uncertainty ellipsoids than for the linearization method, but it encloses few points after 3.5 second. The re-sampling in the unscented
method with re-sampling has the effect of enlarging the uncertainty ellipsoids. Thus,
the size of the uncertainty ellipsoids increase rapidly, and the uncertainty
ellipsoid contains more sample points than the other two methods.

The mean numbers of sample points contained in the computed uncertainty ellipsoid are
10.7\%, 26.3\%, and 73.55\%, respectively. The computation times are 4.97, 43.50, and 45.56
seconds, respectively.

\section{Global Features of the Attitude Flow}

The numerical results presented in the previous section demonstrate
the difficulty in obtaining accurate global bounds on attitude
solutions that are initialized in an uncertainty ellipsoid.    It is
claimed that the source of this difficulty is the nonlinear attitude
flow of the 3D pendulum, especially the fact that the flow can
exhibit chaos and extreme sensitivity to initial conditions.    A
conceptual description of certain global features of the attitude
flow that help to explain this effect is now provided.

As described in~\cite{SIADS07}, the global dynamics of the 3D pendulum are
complicated. There is a 1D hanging equilibrium submanifold of the 3D 
configuration manifold, consisting of hanging equilibria that differ by 
a rotation about the vertical. There is also a 1D inverted equilibrium 
submanifold consisting of inverted equilibria that differ by a rotation 
about the vertical.   Each hanging equilibrium is stable in the sense of 
Liapunov.    Each inverted equilibrium is unstable, with a 2D stable manifold, a 2D
unstable manifold, and a 2D center manifold.   Let $M$ denote the union
of all the 2D stable manifolds corresponding to inverted equlibria. This 3D set 
$M$ plays an important role in understanding the global dynamics of the 
attitude flow.

Every trajectory in $M$ converges to the inverted equilibrium
manifold.   Although the set $M$ has zero measure, its existence
influences the dynamics of the 3D pendulum attitude flow near $M$.
Since $M$ is constructed as the union of the stable manifolds of
unstable equilibria, trajectories near $M$ remain near $M$ for an
extended period of time. In
particular, the closer a trajectory is to $M$ the longer it remains
near $M$.   In fact, there are trajectories that remain close to $M$ for
arbitrarily long time periods.  This property is due to the saddle
character of each inverted equilibrium.

It should be mentioned that it is difficult to determine exactly the
set $M$.   One can make use of linear attitude
equations near an inverted equilibrium to approximate the tangent
space to the stable manifold of that equilibrium. However, this
provides only local information about $M$; the non-local properties of
the set $M$ are not understood. In practice, to accurately compute the global structure of a stable manifold, one relies on either (i) extremely high-order Taylor approximations of the nonlinear stable manifold for a neighborhood of the equilibrium, which is used to obtain sample points on the stable manifold that are then propagated backwards in time in order to compute the global structure of the stable manifold~\cite{GoKoLoMaMaRo2004}, or (ii) set-oriented techniques based on representing the nonlinear flow map for short times as a Markov process~\cite{DeFrJu2001}.

This argument demonstrates that the set $M$ has a strong influence on
the 3D pendulum dynamics near $M$, with high shearing and thus high
sensitivity of the attitude flow near $M$.   This is one of the mechanisms leading to the complex nonlinear dynamics of the 3D pendulum and makes it
impossible to efficiently compute accurate global bounds on attitude solutions
that are initialized in an uncertainty ellipsoid.

\section{Conclusions}

The Lie group variational integrator is known to provide accurate long-term solutions of the rigid body equations in the presence of an external potential for a given initial attitude and angular velocity; these computed solutions exactly conserve the theoretical conservation properties, namely the symplectic structure, and the angular momentum component about the vertical axis in the case for the 3D pendulum.  In addition, it exhibits very good energy behavior, with only a very small bounded energy oscillation, for exponentially long times. Furthermore, the Lie group variational integrator is also known to exactly conserve orthogonality of the computed attitude as a rotation matrix.

Consequently, inaccuracies in the approximate ellipsoidal bounds computed according to the three computational methods introduced do not arise from computational difficulties with the Lie group variational integrator.    Rather, the inaccuracies in the approximate ellipsoidal bounds arise from the fact that the attitude flow dynamics are highly nonlinear, with regions wherein the dynamics cannot be adequately approximated by linear dynamics.

With these qualifications, it is clear that the unscented method with resampling is the most accurate of the three proposed methods in propagating the uncertainty.   This method can provide a basis for analysis of control and estimation problems for attitude systems such as the 3D pendulum.   For example, digital control operates open loop between sample times; the analysis in this paper demonstrates the importance of the choice of inter-sampling time in obtaining accurate propagation of the attitude flow between sample times. As another example, attitude estimation operates open loop between measurement times; the analysis in this paper demonstrates the importance of the choice of inter-measurement time in obtaining accurate propagation of the attitude flow between measurement times.

The bottom line demonstrated by the development in this paper is that guaranteed global bounds for attitude dynamics defined by the 3D pendulum, or indeed for any attitude dynamics with a nontrivial potential, are not achievable.    That is, there is no universal approach to global and uniform approximation of  the attitude flow dynamics.

\bibliography{est}
\bibliographystyle{IEEEtran}

\clearpage\newpage
\renewcommand{\thesubfigure}{}

\renewcommand{\thesubfigure}{}

\end{document}